\documentstyle[12pt]{article}
\input amssym.def
\input amssym.tex

\def\emptyset{\varnothing}
\def\ol{\overline}

\def\ero{{\smallsetminus}}

\def\cub{{\rm CUB}}

\def\dom{{\rm dom}}

\def\alkus{\mathop{\vartriangleleft}}
\def\baire{{}^\kappa \kappa }
\def\bairen{{}^\omega \omega }
\def\cantor{{}^\kappa 2}
\def\cantorn{{}^\omega 2}

\def\Seq{{}^{<\kappa }\kappa }

\def\F{{\cal F}}

\def\Q{{\cal Q}}

\def\T{{\cal T}}

\def\proves{\vdash}

\newtheorem{theorem}{Theorem}[section]
\newtheorem{lemma}[theorem]{Lemma}

\newtheorem{definition}[theorem]{Definition}
\newenvironment{defn}{\begin{definition}\rm }{\end{definition}}
\newtheorem{hypoth}[theorem]{Assumptions}
\newenvironment{hyp}{\begin{hypoth}\rm }{\end{hypoth}}

\newtheorem{rem}[theorem]{Remark}
\newtheorem{ex}[theorem]{Example}
\newtheorem{prob}[theorem]{Problem}

\newenvironment{remark}
               {\begin{rem}\rm}{\end{rem}}

\def\tod{\par\noindent{\it Proof. }}
\def\mot{\nobreak{$\square$}\medskip}

\newenvironment{proof}
               {\tod}{\mot}

\def\theenumi{\arabic{enumi}}

\def\theenumii{\alph{enumii}}
\def\p@enumii{\theenumi}

\def\theenumiii{\Alph{enumiii}}
\def\p@enumiii{\theenumi(\theenumii)}

\def\p@enumiv{\p@enumiii\theenumiii}
\author{Aapo Halko\thanks{ Partially supported by grant \#1011049 of the
Academy of Finland}\and 
Saharon Shelah\thanks{ Publication number 662}}
\date{\today\footnote{1991 Mathematics Subject Classification: 03E15,
04A15.}}
\title{On strong measure zero subsets of ${}^\kappa 2$
}

\begin{document}
\maketitle


\begin{abstract}
This paper answers  three questions posed in \cite{Halk96}.
\begin{description}
  \item[Theorem \ref{thm-additivity}] The family of strong measure zero
subsets of
${}^{\omega _1}2$ is
$2^{{\aleph _1}}$-additive under GMA and CH.
  \item[Theorem \ref{thm-notgbc}] The generalized Borel
conjecture is false in ${}^{\omega _1}2$ assuming ZFC+CH.
  \item[Theorem \ref{thm-bp}] The family of subsets of ${}^{\omega _1}2$ with the
property of Baire is not closed under the Souslin operation.
  \end{description}
\end{abstract}

\section{Introduction}

We study the generalized Cantor space 
$\cantor$ and the generalized Baire space
$\baire$ for an uncountable cardinal $\kappa$ as analogues of the classical
Cantor and Baire spaces. We equip $\baire$ with the
topology where a basic neighborhood of
a point $\eta $ is the set 
$$
\{\nu\in\baire:(\forall j<i)(\nu(j)=\eta(j))\},
$$
where $i<\kappa$. A systematic study of measure and category
in these spaces was started in \cite{Halk96}.
In this paper we answer some problems posed in \cite{Halk96}.

There are natural generalizations of the concepts of 
meager and strong measure zero sets from the
space $\bairen$ to the space $\baire$. Many results and their proofs concerning
these concepts,  e.g. the Baire Categoricity Theorem, are just
straightforward generalizations of the corresponding results of
$\bairen$. It was proved in
\cite{Shel78} that, assuming the Generalized Martin's Axiom
GMA  of \cite{Shel78},  the family of 
meager subsets of ${}^{\omega_1}2$ is closed
under unions of length $<2^{\aleph _1}$. 
In Section \ref{gma} we prove the same additivity result
for the family of strong measure zero sets of $\cantor$.

The generalized Borel conjecture for ${}^\kappa 2$, which we will call by $GBC(\kappa )$, states that every
strong measure zero subset of $\cantor$ has the cardinality at most $\kappa $. The
consistency of the Borel Conjecture for the space
$\cantorn$, i.e. GBC($\omega $), was shown by Laver in \cite{Lave77}. However, in
Section \ref{gbc} we show that GBC($\kappa $)
fails  assuming $\kappa =\kappa ^{<\kappa }=\mu ^+>\aleph _0$. It is an open problem
whether the statements ``$\kappa $ strongly inaccessible + GBC($\kappa $)" or
``$\kappa $ the first (strongly) inaccessible + GBC($\kappa $)" are consistent.

In the final section we show that the property of Baire 
is not preserved by the generalized Souslin operation
\[\bigcup_{f\in\baire}\bigcap_{i<\kappa}A_{f{\restriction}i}.\]
We show this by pointing out that the set
$\cub$ of characteristic functions of closed
unbounded sets of $\kappa $ lacks the property of Baire and yet is obtained from
open sets by this Souslin operation.

We thank Jouko V\"a\"an\"anen for reading this paper and suggesting many
improvements.

Our set theoretical notation is standard, see \cite{Jech78}.
Ordinals are denoted by $\alpha $, $\beta $, $\epsilon $, $\xi $, $i$, $j$; cardinals by $\kappa $, $\mu $ and
sequences by $\eta $, $\nu $. Length of a sequence $\eta $ is denoted by $\ell(\eta )$.  We
denote
$[\alpha ,\beta )=\{ i\mid \alpha \leq i<\beta \} $.  If
$\eta $ and
$\nu $ are sequences, then
$\eta \alkus \nu $ means that
$\eta $ is an initial segment of $\nu $. For a cardinal $\kappa $ and a set $A$ we denote
$[A]^\kappa =\{ B\subseteq A:|B|=\kappa \} $ and $[A]^{\leq \kappa }=\{ B\subseteq A:|B|\leq \kappa \} $.

\section{Strong measure zero sets}\label{gma}

\begin{hyp}\label{hyb}
Assume that $\kappa $ is uncountable. Let $\T\subseteq {}^{<\kappa }\kappa $ be a normal tree with
 $\kappa $ levels. Let $\T_i$ be the $i$-th level of $\T$ and 
$\T_\kappa =\lim_\kappa (\T)$. Assume that
 $$
i<j\leq \kappa  \Rightarrow  (\forall \eta \in \T_i)(\exists \nu \in \T_j)(\eta \alkus \nu ).
 $$ We also assume that
$|i|\leq |\T_i|\leq \kappa $ for each $i<\kappa $ and $|\T_\kappa |>\kappa $.
Let
$F_i:\T_i\to |\T_i|$ be one to one. We denote $\ol F=\langle F_i: i<\kappa \rangle $ and 
$\ol F\circ \eta =\langle F_i(\eta {\restriction}i):i<\kappa \rangle $ for each $\eta \in \T_\kappa $. 
\end{hyp}

\begin{remark}
 If $\kappa =\kappa ^{<\kappa }$, then $\T={}^{<\kappa }2$ and
$\T={}^{<\kappa }\kappa $ satisfy \ref{hyb}. So, in particular, \ref{hyb} is true for
$\T={}^{<\omega _1}\omega _1$ under CH and for $\T={}^{<\kappa }\kappa $ where $\kappa $ is strongly
inaccessible.  
\end{remark}

We introduce some notation. If $\nu \in \T$ then $[\nu ]=\{ \eta \in \T_\kappa :\nu \alkus \eta \} $. For
$X\subseteq \kappa $ and
$f,g\in {}^X\kappa $,
$$
f<^*_\kappa g\iff |\{ i\in X:f(i)\geq g(i)\} |<\kappa .
$$

\begin{defn}
$A\subseteq \T_\kappa $ has {\em strong measure zero\/}, if for every $X\in [\kappa ]^\kappa $ we can find
$\langle f_\xi :\xi \in X\rangle $, $f_\xi \in \T_\xi $ such that 
$$
A\subseteq \bigcup _{\xi \in X}[f_\xi ].
$$
\end{defn}

Next we give two characterizations of strong measure zero sets which we shall
use in the proofs  of the theorems in this and next sections. 

\begin{lemma}\label{lemma-one}
 The following are equivalent for $A\subseteq \T_\kappa $
        \begin{enumerate}
              \item[(a)] $A$ has strong measure zero
              \item[(b)]\label{bee} if $\langle \alpha _i:i<\kappa \rangle $ is strictly increasing
continuous sequence of ordinals $<\kappa $ then we can find 
$$Y_i\in [\T_{\alpha _{i+1}}]^{\leq |\alpha _i|}$$
such that
$$(\forall \eta \in A)(\exists ^\kappa i)(\eta {\restriction}\alpha _{i+1}\in Y_i).$$
 \end{enumerate}
 \end{lemma}
\begin{proof}
(a) implies (b). Let $\langle \alpha _i:i<\kappa \rangle $ be strictly increasing continuous
sequence. For each $i<\kappa $ apply (a) to 
$$
X_i=\{ \alpha _{j+1}:j\geq i\} 
$$
getting $\langle f_{i,\alpha _{j+1}}\in \T_{\alpha _{j+1}}:j\geq i\rangle $. Let 
$$
Y_i=\{ f_{\epsilon ,\alpha _{i+1}}:\epsilon \leq i\} .
$$
Now $|Y_i|\leq |i|\leq |\alpha _i|$ and if $\eta \in A$ then for any $i<\kappa $ there is $j\geq i$ such that  $\eta {\restriction}\alpha _{j+1}=f_{i,\alpha _{j+1}}\in Y_j$.

(b) implies (a). Let $X\in [\kappa ]^\kappa $. Choose by induction on $i<\kappa $, $\gamma _i<\kappa $
such that if $i$ is limit then $\gamma _i=\cup \{ \gamma _j:j<i\} $, and if $i=j+1$ then choose
$\gamma _i>\gamma _j$ such that the set $X_j=[\gamma _j,\gamma _i)\cap X$ has cardinality $|\gamma _j|$. Apply
clause (b) to $\langle \gamma _i:i<\kappa \rangle $: let 
$$
\langle Y_i\in [\T_{\gamma _{i+1}}]^{\leq |\gamma _i|}:i<\kappa \rangle 
$$
be as guaranteed by clause (b). So $|Y_i|\leq |X_i|$ and we let $h_i:Y_i\to X_i$ be
one to one. Let $\langle f_\xi :\xi \in X\rangle $, $f_\xi \in \T_\xi $, be such that if $\xi =h_i(g)$
for
$g\in Y_i$ then $f_\xi =g{\restriction}\xi $. As $[g]\subseteq [f_\xi ]$ we are done.
\end{proof}
  
  \begin{lemma}\label{lemma-two} 
  If $\kappa =\mu ^+$ and $|\T_i|=\kappa $ for $i<\kappa $ large
enough then the following are equivalent for $A\subseteq \T_\kappa $

  \begin{enumerate}
              \item[(a)] $A$ has strong measure zero
  
  \item[(b$'$)] like \ref{lemma-one}(b), but $$Y_i\in [\T_{\alpha _{i+1}}]^{\leq \mu }.$$
  \item[(c)] for every $X\in [\kappa ]^\kappa $, there is $f\in {}^X\kappa $ such that
  $$\neg(f<^*_\kappa (\ol F\circ \eta ){\restriction}X)$$ for each $\eta \in A$.
 \end{enumerate}
\end{lemma}

\begin{proof}
Under the assumptions, \ref{lemma-one}(b) is clearly equivalent to
\ref{lemma-two}(b$'$).

(b$'$) implies (c). Let $X\in [\kappa ]^\kappa $.
We may assume that 
%
if $\alpha \in [\min X,\kappa )$ then $|\T_\alpha |=\kappa $. Let the
closure of $X\cup \{ 0\} $ be enumerated in $\{ \alpha _i:i<\kappa \} $ where $\alpha _i$ are increasing
with $i$. Apply clause (b$'$) and get
$\langle Y_i:i<\kappa \rangle $, $Y_i\in [\T_{\alpha _{i+1}}]^{\leq \mu }$. Choose $f\in {}^X\kappa $ such that
$$
f(\alpha _{i+1})=\min\{ \gamma <\kappa :F_{\alpha _{i+1}}(\eta )<\gamma \hbox{ for every }\eta \in Y_i\} .
$$ 
 Now let $\eta \in A$. Then $H=\{ i<\kappa :\eta {\restriction}\alpha _{i+1}\in Y_i\} $ has cardinality $\kappa $ and 
$F_{\alpha _{i+1}}(\eta {\restriction}\alpha _{i+1})<f(\alpha _{i+1})$ for each $i\in H$.  This means
$\neg(f<^*_\kappa (\ol F\circ \eta ){\restriction}X)$.

(c) implies (b$'$). Let $\langle \alpha _i:i<\kappa \rangle $ be strictly increasing continuous sequence of
ordinals $<\kappa $. We should find $\langle Y_i:i<\kappa \rangle $ as in clause (b$'$). Apply clause~(c)
for
$X=\{ \alpha _{i+1}:i<\kappa \} $ and get $f\in {}^X\kappa $. Let 
$$
Y_i=\{ \eta \in \T_{\alpha _{i+1}}:F_{\alpha _{i+1}}(\eta )\leq f(\alpha _{i+1})\} .
$$
Let $\eta \in A$. Then $H=\{ i<\kappa :F_{\alpha _{i+1}}(\eta {\restriction}\alpha _{i+1})\leq f(\alpha _{i+1})\} $ has cardinality
$\kappa $ and $\eta {\restriction}\alpha _{i+1}\in Y_i$ for all $i\in H$.
\end{proof}

 A family $\F\subseteq {}^\kappa \kappa $ is {\em bounded\/}, if there is $g\in {}^\kappa \kappa $
such that $f<_\kappa ^*g$ for all $f\in \F$.  A family $\F\subseteq {}^\kappa \kappa $ is {\em
dominating\/}, if for each
$g\in {}^\kappa \kappa $ there is $f\in \F$ such that $g<_\kappa ^*f$.
Condition (c) of Lemma \ref{lemma-two} can be rephrased as follows: 
For each $X\in [\kappa ]^\kappa $ the family $\{ (\ol F\circ \eta ){\restriction}X:\eta \in A\} $ is not
dominating. Let ${\bf d}$ be the size of the smallest dominating family and
let ${\bf b}$ be the size of the smallest unbounded family. Clearly $\kappa <{\bf
b}\leq {\bf d}\leq 2^\kappa $.

It is possible to formulate a version of GMA($\kappa $) for arbitrary $\kappa $ with
$\kappa =\kappa ^{<\kappa }$ and prove its relative consistency. See
\cite{Shel78} 1.10 on page 302. 

\begin{lemma}[\cite{Roth48}]\label{gma-fact}
Assume $\kappa =\kappa ^{<\kappa }$ and GMA($\kappa $). Then ${\bf b}=2^\kappa $.
\end{lemma}

We are ready to prove the main result of this section.

\begin{theorem}\label{thm-additivity}
The ideal of strong measure zero sets of ${}^{\kappa }2$ is $2^{\kappa }$-additive 
under $\kappa =\kappa ^{<\kappa }$ and GMA($\kappa $).
\end{theorem}

\begin{proof}
Assume that  $\langle A_\xi :\xi <\gamma \rangle $, $\gamma <2^\kappa $, is a sequence of sets with strong
measure zero. Let $A=\bigcup _{\xi <\gamma }A_\xi $. We prove that $A$ has strong measure zero.
Let
$X\in [\kappa ]^\kappa $. Using (c) of Lemma \ref{lemma-two} for each $\xi <\gamma $ we find
$f_\xi \in {}^X\kappa $ such that
$$
\neg (f_\xi <_\kappa ^*(\ol F\circ \eta ){\restriction}X)
$$
 for all $\eta \in A_\xi $. By Lemma \ref{gma-fact} the set $\{ f_\xi :\xi <\gamma \} $ is bounded.
 Hence there is $f\in {}^X\kappa $ such that 
$$
f_\xi <_\kappa ^*f
$$
for all $\xi <\gamma $. But then 
$$
\neg (f<_\kappa ^*(\ol F\circ \eta ){\restriction}X)
$$
for all $\eta \in A$. Hence $A$ is a strong measure zero set by Lemma
\ref{lemma-two}(c).
\end{proof}

\begin{remark}
Let $\F$ be a dominating family of size ${\bf d}$. Let $X\in [\kappa ]^\kappa $ be such that
$X$ contains no limit ordinals. For each $f\in \F$ we can find $\eta _f\in \baire $ such that
$f<^*_\kappa (\ol F\circ\eta _f){\restriction}X$. Now the set $A=\{ \eta _f:f\in
\F\} $ does not have strong measure zero by Lemma \ref{lemma-two}. Hence the ideal
of strong measure zero sets is not
${\bf d}^+$-additive. So consistently, $\kappa =\kappa ^{<\kappa }$, the ideal is not
$\kappa ^{++}$-additive and $\kappa ^{++}\leq 2^\kappa $.
\end{remark}

\section{The generalized Borel conjecture}\label{gbc}

Let the {\em Generalized Borel Conjecture for $\T_\kappa $\/} 
be the statement that every strong measure zero subset of $\T_\kappa $ has
cardinality at the most $\kappa $. Let GBC($\kappa $) be the generalized Borel
conjecture for $\cantor$ and let  GBC be GBC($\aleph _1$).

\begin{theorem}\label{thm-notgbc}
$ZFC + CH \proves \neg GBC$.
\end{theorem}

This theorem follows from the following more general lemma. 

\begin{lemma}
If $\kappa =\kappa ^{<\kappa }=\mu ^+$, $|\T_i|=\kappa $ for $i<\kappa $ large enough and $\T$ is closed under
increasing sequences of length $<\kappa $ then there is an $A\in [\T_\kappa ]^{\kappa ^+}$ of strong
measure zero.
\end{lemma}

\begin{proof}
We consider two cases, according to the size of cardinal number {\bf d}.

{\bf Case 1:} ${\bf d}>\kappa ^+$.
Let $A\subseteq \T_\kappa $ be any set of cardinality $\kappa ^+$. We shall prove it has strong
measure zero. Let $X\in [\kappa ]^\kappa $. The set $\{ (\ol F\circ \eta ){\restriction}X:\eta \in A\} $ is not
dominating in
$({}^X\kappa ,<^*_\kappa )$. Hence there is
$f\in {}^X\kappa $ such that $\neg(f<^*_\kappa \ol (F\circ \eta ){\restriction}X)$ for every $\eta \in A$. But then
$A$ has strong measure zero by  clause (c) of Lemma \ref{lemma-two}.

{\bf Case 2:} ${\bf d}=\kappa ^+$.
Let $\langle g_\epsilon ^*:\epsilon <\kappa ^+\rangle $ be dominating. We may assume that each $g_\epsilon ^*$ is increasing
and if $\epsilon <\zeta $ then $g_\epsilon ^*<^*_\kappa g_\zeta ^*$. Let 
$$
C_\epsilon ^*=\{ \delta <\kappa :\delta \hbox{ limit}\land\forall i(i<\delta  \iff  g_\epsilon ^*(i)<\delta )\} .
$$ Let
 $C_\epsilon ^*=\{ \alpha _{\epsilon ,i}^*:i<\kappa \} $ where $\alpha _{\epsilon ,i}^*$ is increasing in $i$. We choose 
$\eta _\epsilon \in \T_\kappa $ and   $\langle Y_{\epsilon ,i}:i<\kappa \rangle $ by induction on $\epsilon <\kappa ^+$ such that
\begin{enumerate}
  \item $\eta _\epsilon \not\in \{ \eta _\zeta :\zeta <\epsilon \} $
  \item $Y_{\epsilon ,i}\in [\T_{\alpha _{\epsilon ,i+1}^*}]^{\leq \mu }$
  \item if $\zeta \leq \epsilon $ then $(\exists ^\kappa i<\kappa )(\eta _\zeta {\restriction}\alpha _{\epsilon ,i+1}^*\in Y_{\epsilon ,i})$
  \item\label{rajaehto} for every $\nu \in \T$ and for every $i<\kappa $ large enough 
  $$
  (\exists \rho )(\nu \alkus \rho  \land \rho \in Y_{\epsilon ,i})
$$ 
  \item\label{v-ehto} if $\zeta >\epsilon $ then $(\exists ^\kappa i<\kappa )(\eta _\zeta {\restriction}\alpha _{\epsilon ,i+1}^*\in Y_{\epsilon ,i})$.
 \end{enumerate}
This can be done as follows:

Choose for each $\nu \in \T$ some $\rho _\nu \in [\nu ]$ and let
$\Q=\{ \rho _\nu :\nu \in \T\} $. Since we assume $\kappa =\kappa ^{<\kappa }$ we can enumerate $\Q$ in
$\{ \rho _i:i<\kappa \} $. For a start, let $\eta _\zeta =\rho _\zeta $ when $\zeta <\kappa $ and
$$
Y_{\zeta ,i}=\{ \eta _j{\restriction}\alpha _{\zeta ,i+1}^*:j\leq i\} 
$$
for $\zeta <\kappa $ and $i<\kappa $.  Conditions (1)--(5) hold so far.
Assume that $\eta _\zeta $ and $Y_{\zeta ,i}$ have been defined for $\kappa \leq \zeta <\epsilon $ and $i<\kappa $. We will
define
$\eta _\epsilon $ and $Y_{\epsilon ,i}$ for $i<\kappa $ as follows.  We will define certain ordinals
$\beta _j^\epsilon $ and restrictions $\eta _\epsilon {\restriction}\beta _j^\epsilon $ by induction on $j<\kappa $ such that
$\langle \beta _j^\epsilon :j<\kappa \rangle $ is a strictly increasing continuous sequence converging to $\kappa $.
Let $\pi _\epsilon :\epsilon \times \kappa \to \kappa $ be a bijection. If $\beta ^\epsilon _j$ and thereby $\eta _\epsilon {\restriction}\beta _j^\epsilon $ are
defined, let $\beta ^\epsilon _{j+1}$ be as follows. If
$\pi _\epsilon (\zeta ,i)=j$ then let $\nu \in \T$ be such that
$$
\eta _\epsilon {\restriction}\beta ^\epsilon _j\alkus \nu  \hbox{ and $\nu $ not compatible with $\eta _\zeta $.}
$$
 By (\ref{rajaehto}) there is $\rho _i^\zeta \in Y_{\zeta ,j_i^\zeta }$ for some $j_i^\zeta <\kappa $ such that
$\nu \alkus \rho _i^\zeta $.  Let $\beta ^\epsilon _{j+1}=\ell(\rho _i^\zeta )=\alpha _{\zeta ,j_i^\zeta +1}^*$ and
$\eta _\epsilon {\restriction}\beta ^\epsilon _{j+1}=\rho _i^\zeta $. Let $\theta _\epsilon :\kappa \to \epsilon $ be a bijection and
$$
Y_{\epsilon ,i}=\{ \rho _j{\restriction}\alpha _{\epsilon ,i+1}^*:j\leq i\} \cup \{ \eta _{\theta _\epsilon (j)}{\restriction}\alpha _{\epsilon ,i+1}^*:j<i\} \cup \{ \eta _\epsilon {\restriction}\alpha _{\epsilon ,i+1}^*\} .
$$
Conditions (1)--(3) hold trivially.  To see
(\ref{rajaehto}), let $\nu \in \T$. Let $i_\nu <\kappa $ be such that $\alpha ^*_{\epsilon ,i_\nu }\geq \ell(\nu )$
and $i_\nu \geq j_\nu $ where $j_\nu $ is such that $\rho _{j_\nu }=\rho _\nu $. Hence 
$$
(\forall i>i_\nu )(\nu \alkus \rho _{j_\nu }{\restriction}\alpha _{\epsilon ,i+1}^*\land \rho _{j_\nu }{\restriction}\alpha _{\epsilon ,i+1}^*\in Y_{\epsilon ,i}).
$$
For (\ref{v-ehto}), let $\zeta >\epsilon $. By construction 
$$
\eta _\zeta {\restriction}\alpha _{\epsilon ,j_i^\epsilon +1}^*\in Y_{\epsilon ,j_i^\epsilon }
$$
for all $i<\kappa $. Hence 
 $$
 \exists ^\kappa i(\eta _\zeta {\restriction}\alpha _{\epsilon ,i+1}^*\in Y_{\epsilon ,i}).
 $$

Let $A=\{ \eta _\epsilon :\epsilon <\kappa ^+\} $. Clearly $|A|=\kappa ^+$. Now we show that $A$ is of
strong measure zero by using clause (b$'$) of Lemma \ref{lemma-two}: Let
$\langle \alpha _i:i<\kappa \rangle $ be a club where
$\alpha _i$ is increasing with
$i$. So for some $\epsilon $, $g_\epsilon ^*$ dominates $i\mapsto \alpha _i$.   Let $i_0$ be
such that $(\forall i>i_0)(\alpha _i<g^*_\epsilon (i))$. If $i$ is such that
$\alpha _{\epsilon ,i}^*>i_0$ then
$\alpha _{\epsilon ,i}^*=\sup\{ \alpha _j:j<\alpha _{\epsilon ,i}^*\} $ because
$i_0<j<\alpha _{\epsilon ,i}^*$ implies $\alpha _j<g^*_\epsilon (j)<\alpha _{\epsilon ,i}^*$. Hence for every
$i$ large enough $\alpha _{\epsilon ,i}^*\in \{ \alpha _j:j<\kappa \} $.
Define
$$
Y_i=\{ \rho {\restriction}\alpha _{i+1}:\rho \in Y_{\epsilon ,j}\hbox{ where $j$ is minimal such that
}\alpha _{\epsilon ,j+1}^*\geq \alpha _{i+1}\} .
$$
$\langle Y_i:i<\kappa \rangle $ is as required: 
Clearly $|Y_i|\leq \mu $. Suppose $\eta _\zeta \in A$. By (3) and (5), there is a strictly
increasing sequence $\langle i_\xi :\xi <\kappa \rangle $ in $\kappa $ such that 
$$
\eta _\zeta {\restriction}\alpha ^*_{\epsilon ,i_\xi +1}\in Y_{\epsilon ,i_\xi }
$$
for all $\xi <\kappa $. Choose $j_\xi $ such that 
$$
\alpha ^*_{\epsilon ,i_\xi }<\alpha _{j_\xi +1}\leq \alpha ^*_{\epsilon ,i_\xi +1}.
$$
Clearly, $\xi <\xi '$ implies $j_\xi <j_{\xi '}$. Now 
$
\eta _\zeta {\restriction}\alpha _{j_\xi +1}\in Y_{j_\xi }$ for all $\xi <\kappa $ and the claim follows.
\end{proof}

\section{The property of Baire}\label{pb}

The topology of $\cantor$ is the one generated by the $[\eta ]$ as basic
neighborhoods. So 
$A\subseteq {}^\kappa 2$ is {\em open\/}, if for every $\eta \in A$ there is $i<\kappa $ such that
$[\eta {\restriction}i]\subseteq A$. 
$A$ is {\em nowhere dense\/}, if for every $\nu \in {}^{<\kappa }2$ there is
$\eta \in {}^{<\kappa }2$ such that $\nu \alkus \eta $ and $A\cap [\eta ]=\emptyset $. $A$ is {\em meager\/}, if
$A=\bigcup _{\xi <\kappa }R_\xi $ where the sets $R_\xi \subseteq {}^\kappa 2$ are nowhere dense. $A$ has {\em the
property of Baire\/}, if there is an open set $O\subseteq {}^\kappa 2$ such that $(O\ero
A)\cup (A\ero O)$ is meager.


Let $$
\cub=\{ \eta \in {}^\kappa 2:\hbox{ for some club $C$ of }\kappa \ (\forall i\in C)(\eta (i)=1)\} .
$$
\begin{lemma}
 There is a system 
$\langle A_\nu :\nu \in \Seq\rangle $ of open sets such that
 $$
 \cub=\bigcup _{f\in \baire}\bigcap _{i<\kappa }A_{f{\restriction}i}
$$
 \end{lemma}

 \begin{proof}
 For $\nu \in \Seq$ let
 $$
 A_\nu =\{ \eta \in \cantor:(\forall i\in \dom(\nu ))(\eta (\nu (i))=1)\} 
 $$ 
if $\nu $ is a strictly increasing continuous sequence  and let $A_\nu $ be empty
otherwise. Let $\eta \in \cub$ and let $\langle \alpha _i:i<\kappa \rangle $ be an increasing enumeration of a
club set such that $\eta (\alpha _i)=1$ for all
$i<\kappa $. Then
$\eta \in A_{\langle \alpha _j:j<i\rangle }$ for all
$i$. Conversely, if
$\eta \in A_{f{\restriction}i}$ for all $i$, then clearly $f$ is strictly increasing and
continuous, hence
$\eta \in \cub$.
 \end{proof}

The above lemma shows that the set \cub\ can be obtained from open sets by
means of an operation which is analogous to the Souslin operation. Thus the
following result shows that the property of Baire is not preserved by this
``Souslin" operation. Recall that in the space $\cantorn$ the property of
Baire is preserved by the ordinary Souslin operation.

\begin{theorem}\label{thm-bp}
Let $\kappa 
>\aleph _0$ be regular. 
Then \cub\ does not have the property of Baire.
\end{theorem}

\begin{proof}
We show that for all  open set $O$,
$(O\ero\cub)\cup (\cub\ero O)$ is not meager.

Suppose first $O$ is empty. We show that $\cub$ is not meager. Let $R_\xi \subseteq {}^\kappa 2$
be nowhere dense for
$\xi <\kappa $. We choose
$\alpha _i$,
$\eta _i$ by induction on $i\leq \kappa $ such that
\begin{enumerate}
  \item $\eta _i\in {}^{\alpha _i}2$
  \item if $j<i$ then $\alpha _j<\alpha _i$ and $\eta _j\alkus \eta _i$
  \item if $i$ is limit then $\alpha _i=\bigcup _{j<i}\alpha _j$ and $\eta _i=\bigcup _{j<i}\eta _j$ 
  \item $\eta _{i+1}(\alpha _i)=1$
  \item $\neg(\exists \rho )(\eta _{i+1}\alkus \rho \land \rho \in R_i)$.
 \end{enumerate}
 Now  $\eta _\kappa \in \cub\ero\bigcup _{\xi <\kappa }R_\xi $, whence $\cub\not=\bigcup _{\xi <\kappa }R_\xi $. 
 
If $O$ is non-empty then we choose $\nu $ such that $[\nu ]\subseteq O$. Then
$O\ero\cub\supseteq [\nu ]\ero\cub$.  Similarly as above we show that $[\nu ]\ero
\cub$ is not meager.
We proceed as above except
$\alpha _0=\ell(\nu )$,
$\eta _0=\nu $ and
\begin{enumerate}
  \item[(4$'$)] $\eta _{i+1}(\alpha _i)=0$.
 \end{enumerate}
 Then $\eta _\kappa \in ([\nu ]\ero\cub)\ero\bigcup _{\xi <\kappa }R_\xi $. 
 \end{proof}

Let us call a subset of $\cantor$ {\em Borel\/} if it is a member of the
smallest algebra of subsets of $\cantor$ containing all open sets and closed
under complements and unions of length $\leq \kappa $. It is proved in \cite{Halk96}
that Borel sets have the property of Baire. Hence $\cub$ is not Borel. This
improves the result in \cite{MeVa93} to the effect that $\cub$ is not $\Pi _3^0$
or $\Sigma _3^0$. 
Assuming $\kappa =\aleph _1=2^{\aleph _0}$, non-Borelness of $\cub$ follows from the stronger
result that $\cub$ and ${\rm NON\hbox{-}STAT}=\{ \eta \in {}^{\omega _1}2:\hbox{for some cub
}C\subseteq \omega _1(\forall i\in C)(\eta (i)=0)\} $ cannot be separated by a Borel set \cite{ShVa}.

Aapo Halko

Department of Mathematics

P.O. BOX 4

FIN-00014 University of Helsinki

Helsinki, Finland
 
{\tt aapo.halko@helsinki.fi}

\vspace{1cm}
Saharon Shelah

Institute of Mathematics

Hebrew University

Jerusalem, Israel

\end{document}